\documentclass[12pt,oneside]{article}
\usepackage{amsmath,amssymb,xypic}
\def\R{{\mathbb R}}

\def\t{\widetilde}
\def\ge{\geqslant}

\begin{document}
%

\begin{center}
\large\bf
On the contact equivalence problem of second order ODEs
which are quadratic with respect to the second order derivative
\\[5pt]
\normalsize
\rm Vadim V. Shurygin, jr.\\[5pt]
{\it Kazan (Volga Region) Federal University, Russia}\\
{\it e-mail:} {\tt vshjr@yandex.ru, vadimjr@kpfu.ru}
\end{center}

{\bf Abstract.}
{\small
In the present paper we establish the necessary and sufficient
conditions for two ordinary differential equations of the form
$y''{}^2+A(x,y,y') y''+B(x,y,y')=0$ to be equivalent
under the action of the pseudogroup of contact transformations.
These conditions are formulated in terms of integrals of some
one-dimensional distributions.


\vskip1cm


\section{Introduction}

 In the present paper we solve the problem of equivalence
of 2nd order ODEs of the form
\begin{equation}
\label{du2}
y''{}^2+A(x,y,y') y''+B(x,y,y')=0
\end{equation}
under the action of the pseudogroup of contact transformations of
the space
$\R^3(x,y,p)$, where $p=y'$.
We use the results of A.\,Tresse and B.\,Kruglikov
who described the algebra of differential invariants of 2nd order
ODEs of the form
\begin{equation}
\label{du1}
y''=f(x,y,y')
\end{equation}
under the action of the pseudogroup of point transformations of
the space
$\R^2(x,y)$.

S.\,Lie showed that any two second order ODEs of the form
(\ref{du1})
are equivalent under the action of contact transformations pseudogroup
of the space $\R^3(x,y,p)$.
The problem of point equivalence of such ODEs was the point of
investigation of many papers.
Lie proved that if the ODE~(\ref{du1}) is point equivalent to
the linear equation $y''=0$, it is necessary cubic with respect
to the first derivative, i.e., is of the form
$$
y''=a_3(x,y)y'^3+a_2(x,y)y'^2+a_1(x,y)y'+a_0(x,y).
$$
The class of such equations is closed under the action of point
transformations pseudogroup. Moreover, he formulated the
sufficient conditions of the ODE to be linearizable
(see~\cite{Lie}). R.\,Liou\-ville~\cite{Lio} found precise
conditions for linearization. A.\,Tresse~\cite{Tr2} found the
complete set of relative differential invariants of
ODEs~(\ref{du1}). B.\,Kruglikov~\cite{Krug} formulated the results
of Tresse in a modern language. He described the algebra of
absolute differential invariants and solved the problem of
equivalence of such ODEs under the action of point transformations
pseudogroup. The class of ODEs, cubic with respect to the first
derivatives was considered by Tresse~\cite{Tr1},  who also found
the  relative invariants of these ODEs.
V.\,Yumaguzhin~\cite{Yuma1, Yuma2} studied the the problem of
point equivalence of such ODEs and the obstructions to
linearization. We also note the book of N.\,Ibragimov~\cite{Ibr},
the paper of  N.\,Ibragimov and F.\,Magri~\cite{IM} and the paper
of O.\,Morozov~\cite{Mor}.


\section{Preliminaries}

Let  $\pi$ be a vector bundle and $G$ be the pseudogroup of
diffeomorphisms acting on $\pi$. The action of $G$ naturally
prolongs to the action on the space  $J^k\pi$ of $k$-jets of the
sections of $\pi$. The function $I\in C^\infty(J^k\pi)$ is called
an {\it (absolute) scalar differential invariant of order $k$} if
it is constant along the orbits of prolonged action of  $G$ on
$J^k\pi$. The set of all absolute differential invariants is an
algebra $\cal A$ with respect to the standard arithmetic
operations.

The function $F\in C^\infty(J^k\pi)$ is called a
{\it relative scalar differential invariant of order $k$}
if for every $g\in G$ one has
$g^* F=\mu(g)\cdot F $, for a smooth function $\mu:G\to
C^\infty(J^k\pi)$ satisfying the conditions
 $\mu(g\cdot h)=
h^*\mu(g)\cdot \mu(h)$, $\mu(e)=1$.
In other words, the equation $F=0$ is invariant under the action
of $G$.
The function  $\mu$ is called a {\it weight}.
Let ${\cal M}=\{\mu\}$ be the space of all weights.
Denote by ${\cal R}^\mu$ the space of relative differential invariants
of weight $\mu$.
Then ${\cal R}=\bigcup\limits_{\mu\in{\cal M}}{\cal R}^\mu$ is a
${\cal M}$-graded module over the algebra of
absolute differential invariants.
An {\it invariant differential operator} is a combination
of total derivatives  $\mathcal{D}_{x^i}$ in base coordinates
($x^1,\dots, x^n)$
of  $\pi$
$$
\Delta=\sum\limits^{n}_{i=1}
A_i
{\cal D}_{x^i} + A,\qquad A_i, A \in C^\infty(J^k\pi),
$$
which commutes with the action of $G$.

Let us formulate the  classification of Tresse and Kruglikov which
we will need in what follows.

Any ODE $y''=f(x,y,y')$ may be represented as a section of the bundle
 $J^0\R^3(x,y,p)=\R^4(x,y,p,q)$, where
$p=y'$, $q=y''$ denote the corresponding coordinates in jet
space.
Let ${\cal D}_x$, ${\cal D}_y$, ${\cal D}_p$ be the total
derivatives in $x$, $y$ and $p$ respectively.
We also denote $\hat{\cal D}_x={\cal D}_x+p\,{\cal D}_y$.
Let $q_{lm}^k=\hat{\cal D}_x^l{\cal D}_y^m{\cal D}_p^k(q)$.

The weights form a two-dimensional lattice
and the space of relative invariants is a direct sum
${\cal R}=\bigoplus_{r,s\in \mathbb{Z}}{\cal R}^{r,s}$ (see~\cite{Krug}).
First relative invariants appear in order 4, they are
\begin{multline}
\label{IH}
I=q^4_{00},\qquad
H=
q^2_{20}-4q^1_{11}+6q_{02}^0+q(2q^3_{10}-3q^2_{01})\,- \\[3pt]
- q^1_{00}(q^2_{10}-4q^1_{01})+q^3_{00}q_{10}^0-3q^2_{00}q_{01}^0+q\cdot q\cdot
q^4_{00}.
\end{multline}
Invariant differential operators are

\begin{equation}
\label{Delta_xyp}
\begin{array}{ll}
\Delta_{ p}= & \mathcal{D}_p + (r-s)\dfrac{q^5_{00}}{5q^4_{00}}:
\mathcal{R}^{r,s}\to \mathcal{R}^{r-1,s+1},
\\[10pt]
\Delta_{ x}= & \mathcal{D}_x + p\mathcal{D}_y +
q\mathcal{D}_p +
r\biggl(3q^1_{00}+2\dfrac{q^5_{00}q+q_{10}^4}{q^4_{00}}\biggr)+
\\[8pt]
&
+\, s \biggl(2q^1_{00}+\dfrac{q^5_{00}q+q_{10}^4}{q^4_{00}}\biggr):
\mathcal{R}^{r,s}\to \mathcal{R}^{r+1,s},
\\[10pt]
\Delta_{ y}= &
\dfrac{q^5_{00}}{5q^4_{00}}\mathcal{D}_x +
\biggl(1+p\dfrac{q^5_{00}}{5q^4_{00}}\biggr)\mathcal{D}_y
+ \biggl(2q^1_{00}+\dfrac{5q_{10}^4 + 6q^5_{00}q}{q^4_{00}}\biggr)\mathcal{D}_p
+
\\[8pt]
&
+\, r\biggl(\dfrac{3q^2_{00}}{8} +\dfrac{q_{01}^4}{4q^4_{00}}
+ \dfrac{19q^1_{00}q^5_{00}}{10q^4_{00}}
+ \dfrac{21(q^5_{00}q+q_{01}^4)q^5_{00}}{20q^4_{00}q^4_{00}} \biggr)
+
\\[8pt]
& +\, s \biggl(\dfrac{q^2_{00}}{4}+ \dfrac{q_{01}^4}{2q^4_{00}} +
\dfrac{3q^1_{00}q^5_{00}}{5q^4_{00}}
+ \dfrac{3(q^5_{00}q+q_{01}^4)q^5_{00}}{10q^4_{00}q^4_{00}}\biggr):
\mathcal{R}^{r,s}\to \mathcal{R}^{r,s+1}.
\end{array}
\end{equation}

We say that the equation $y''=f(x,y,y')$ is {\it generic}, if $IH\ne0$ for it.

\smallskip
{\bf Theorem 1.~\cite{Krug}}
{\it
The space of relative differential invariants
of a generic equation
$y''=f(x,y,y')$ is generated by an invariant $H$ and the invariant
differential operators
$\Delta_x$,
$\Delta_y$, $\Delta_p$.}
\smallskip

We also formulate  the equivalence theorem for the generic ODEs.

The functions
$$
H_{10}=\Delta_x H,\quad
H_{01}=\Delta_y H, \quad
K=\Delta_p H
$$
form the basis of invariants of order 5.
There are 11 basic invariants in order 6, we need 8 out of them,
namely
$$
\begin{array}{l}
H_{20}=\Delta_x^2 H, \quad
H_{11}=\Delta_x\Delta_y H,\quad
H_{02}=\Delta_y^2 H,\quad
K_{10}=\Delta_x K,\quad
K_{01}=\Delta_y K, \\[5pt]
\Omega^6=q^6_{00}-\dfrac65\cdot\dfrac{q^5_{00}\cdot
q^5_{00}}{q^4_{00}},
\quad
\Omega^5_{10}=\dfrac{ 5I}{24H}([\Delta_p,\Delta_x]H -
\Delta_yH),\quad
\Omega^4_{20}=\Delta_p^2H -
\dfrac{\Omega^6}{5I}H.
\end{array}
$$
We also denote
$$
J_1 = I^{-1/8}H^{3/8}\in{\cal R}^{1,0},~~
 J_2 = I^{1/4}H^{1/4}\in{\cal R}^{0,1}.
$$
Using $J_1$ and $J_2$ construct the following invariants of pure
order 5 and 6:
$$
\begin{array}{l}
\bar H_{10}=H_{10}/(J_1^3J_2),~~
\bar H_{01}=H_{01}/(J_1^2J_2^2),~~
\bar K=K/(J_1J_2^2),~~\\[3pt]
\bar H_{20}=H_{20}/(J_1^4J_2),~~
\bar H_{11}=H_{11}/(J_1^3J_2^2),~~
\bar H_{02}=H_{02}/(J_1^2J_2^3),~~\\[3pt]
\bar K_{10}=K_{10}/(J_1^2J_2^2),~~
\bar K_{01}=K_{01}/(J_1J_2^3),~~\\[3pt]
\bar \Omega^6=\Omega^6/(J_1^{-4}J_2^5),~~
\bar \Omega^5_{10}=\Omega^5_{10}/(J_1^{-2}J_2^4),~~
\bar \Omega^4_{20}=\Omega^4_{20}/(J_2^3).
\end{array}
$$

Any second order ODE $\cal E$ may be considered as a section
${\mathfrak s}_{\cal E}$ of a bundle $J^0\R^3(x,y,p)$.
Hence, any differential invariant $I$  of order $k$
may be restricted to the equation $\cal E$
via pull-back of the $k$th prolongation of the section ${\mathfrak s}_{\cal E}$:
$$
I^{\cal E}:=({\mathfrak s}_{\cal E}^{(k)})^*(I).
$$

Let  $\cal E$ be a generic second order ODE such that the
functions
$\bar H_{10}^{\cal E}$, $\bar H_{01}^{\cal E}$, $\bar K^{\cal E}$
are local coordinates on $\R^3(x,y,p)$.
The other differential invariants can be expressed as functions of
these three:
$$
\bar H_{ij}^{\cal E}=\Phi_{ij}^{\cal E}(\bar H_{10}^{\cal E},
\bar H_{01}^{\cal E}, \bar K^{\cal E}),
~~
\bar K_{ij}^{\cal E}=\Psi_{ij}^{\cal E}(\bar H_{10}^{\cal E},
\bar H_{01}^{\cal E}, \bar K^{\cal E}),
~~
\bar \Omega_{ij}^{k\cal E}=\Upsilon_{ij}^{k\cal E}
(\bar H_{10}^{\cal E}, \bar H_{01}^{\cal E}, \bar K^{\cal
E}).
$$

\smallskip
{\bf Theorem 2.~\cite{Krug}}
{\it
Two generic second order ODEs ${\cal E}_1$ and  ${\cal E}_2$ are
equivalent under the action of the pseudogroup of point
transformations if and only if the functions
$$
\Phi_{20}^{\cal E}, \Phi_{11}^{\cal E}, \Phi_{02}^{\cal E},
\Psi_{10}^{\cal E}, \Psi_{01}^{\cal E},
\Upsilon^{6\cal E}, \Upsilon_{10}^{5\cal E},
\Upsilon_{20}^{4\cal E}
$$
coincide}.
\smallskip

The criteria for equivalence of non-generic ODEs can be found in~\cite{Krug}.


\section{Various formulations of the problem}

Let $M$ be a 3-dimensional contact manifold with the contact
distribution $\mathcal{ C}$. Let $(x,y,p)$ denote the canonical
coordinates in $M$. Then
 $\mathcal{C}$ is determined by a contact 1-form
$\omega=dy- p\, dx$.

We will say the the ODE (\ref{du2})
is of  {\it hyperbolic type} if
$D=A^2-4B>0$, and that it is of
{\it of elliptic type} if
$D<0$.
Since we are interested in the real case, we restrict ourselves to
the ODEs of hyperbolic type.
Let
\begin{equation}
\label{l1l2}
y''=\lambda_1(x,y,y') \qquad\hbox{and}\qquad y''=\lambda_2(x,y,y')
\end{equation}
be the roots of this ODE.
Note that  $D$ is a relative contact invariant
on the set of ODEs (\ref{du2}).

Geometrically, each of the ODEs $y''=\lambda_i(x,y,y')$
determines a 1-dimensional distribution
${\cal F}_i$ in
$\R^3(x,y,p)=J^0\R^2(x,y)$.
This
distribution lies in $\cal C$ and is spanned by a vector field
\begin{equation}
\label{Xi}
X_i=\dfrac{\partial}{\partial x} + p
\dfrac{\partial}{\partial y } +
\lambda_i(x,y,p) \dfrac{\partial}{\partial p}.
\end{equation}
The Lie bracket
$[X_1, X_2]$ does not belong to $\cal C$.

{\bf Definition.}
We will say that two 1-dimensional distributions ${\cal F}_1$
and ${\cal F}_2$
in a 3-dimensional manifold $M$
are  {\it in general position} if

1) vector fields $X$, $Y$ which define these distributions, are
linear independent at any point;

2) their Lie bracket $[X,Y]$ is linear independent with $X$, $Y$.

Note that the 2-dimensional distribution ${\cal F}_1\oplus {\cal F}_2$
determines the contact structure on $M$.

\smallskip
{\bf Theorem 3.}
{\it The ODE {\rm (\ref{du2})} of hyperbolic type
is equivalent to a pair of 1-dimensional distributions
in general position  in $\R^3(x,y,p)$.}

{\bf Proof.}
The necessity is obvious.

Let ${\cal F}_1$ and ${\cal F}_2$ be two 1-dimensional distributions
in general position in $\R^3(x,y,p)$ and let $X$, $Y$
be vector fields determining these distributions.
Denote ${\cal F}={\cal F}_1\oplus {\cal F}_2$.
Let us prove that the condition 2) of the definition does not depend on a choice
of
the basis in $\cal F$.
For any functions $f$ and $g$
we have
$$
[fX, gY]=fg\cdot[X,Y]+ f\cdot Xg\cdot Y -g\cdot Yf \cdot X =
fg \cdot[X,Y] ~{\rm mod}~{\cal F}.
$$
It follows that for another basis $X'=\alpha X+\beta Y$, $Y'=\gamma X + \delta
Y$ one has
$$
[X', Y']=
(\alpha\delta-\beta\gamma) \cdot[X,Y] \ne 0 ~{\rm mod}~{\cal F},
$$
since $[X,Y]\ne 0~{\rm mod}~{\cal F}$,
$\alpha\delta-\beta\gamma\ne 0$.

Let $a$, $b$ be independent integrals of $X$ and $f$, $g$
be independent integrals of $Y$.
Since  $X$ and $Y$ are linear independent, any three of these
functions are independent and the fourth one can be expressed
in terms of these three. Let, for example, $g=h(a,b,f)$.
Choose the functions  $a$, $b$, $f$ as local coordinates.
In these coordinates the fields $X$ and $Y$ (up to the factors) have the form
$$
X=\dfrac{\partial}{\partial f},
\quad
Y=h_{b}\dfrac{\partial}{\partial a} -
h_{a}\dfrac{\partial}{\partial b}.
$$
Note that  $h_a\ne0$,
since $b$, $f$, $g$ are independent.
Similarly,
$h_b\ne0$.
Denote $c=-{h_{a}}/{h_{b}}$.
Since
$$
[X, Y]=h_{bf}\dfrac{\partial}{\partial a}
- h_{af}\dfrac{\partial}{\partial b}
$$
does not belong to $\cal F$, one has
$$
h_{a}h_{bf} - h_{b}h_{af}\ne 0.
$$
This is equivalent to the inequality
$$
\dfrac{\partial}{\partial f}\left(\dfrac{h_{a}}{h_{b}}\right)\ne0,
$$
which means that
one can choose
$a$, $b$, $c$ as coordinates on $\R^3$.
In these coordinates  $X$, $Y$
will have the form (up to the factors)
$$
X=\dfrac{\partial}{\partial c},~~
Y=\dfrac{\partial}{\partial a} +
c\dfrac{\partial}{\partial b}- \dfrac{{\rm Flex} (h)}{h_{b}^3}
\dfrac{\partial}{\partial c},
$$
where we denoted
$$
{\rm Flex} (h)= h_{aa}h_{b}^2 -2h_{a}h_{b}h_{ab} + h_{bb}h_{a}^2.
$$
The distribution
$\cal F$ is determined by the
1-form $\omega=db- c\, da$ and defines the contact structure on $\R^3$.

Note that if a vector field $Z$ lies in  $\cal F$ and is linear independent with $X$,
then up to a factor it is of the form~(\ref{Xi}).
Consider the contact transformation which satisfies the following
conditions:
1) it preserves some point in $\R^3$, 2) its differential is
sufficiently close to identity and 3) the direction of the field $X$
is not an eigendirection.
Such a transformation maps
$X$ and $Y$ to the fields which are not proportional to $X$.
These vector fields up to the factors will have the form
(\ref{Xi}). Thus, they determine two equations
(\ref{l1l2}).
$\Box$
\smallskip

Let us also show that the ODE (\ref{du2}) can be treated as
two pairs of functions on $M$, satisfying some additional conditions.

The equation $y''=\lambda_1(x,y,y')$ locally has
a pair of integrals  $(a,b)$  (which are also integrals of $X_1$) such that
$da\wedge db\ne 0$.
Since $X_1$ lies in the kernel of the forms $\omega$, $da$ and
$db$, one has
$\omega=\alpha\, da + \beta \, db$ for some
functions $\alpha$, $\beta$.
Thus, ODE (\ref{du2}) of hyperbolic type, determines
two pairs of functions
$(a,b)$ and $(f,g)$ such that any three of the are independent.
It follows from the equality $da(X_1)=X_1a$ that
1-forms $(X_1g)\, da- (X_1f)\, db$ и
$(X_2b)\, df- (X_2a)\, dg$
are proportional to $\omega$.

{\bf Definition.}
Let  $(a,b)$ and $(f,g)$ be two pairs of functions on $M$.
We will say that they are
{\it in general position}, if:

1) any three of them are independent;

2) the 1-dimensional distributions determined by $X$ and $Y$
are in general position.

\smallskip
{\bf Proposition 1.}
{\it Let two pairs of functions $(a,b)$ and $(f,g)$ on $M$
be in general position. Let $a'=a'(a,b)$, $b'=b'(a,b)$ be independent functions,
and $f'=f'(f,g)$, $g'=g'(f,g)$ be independent functions.
Then the pairs $(a',b')$ and $(f',g')$ also are in general position.}

{\bf Proof.}
Follows from the fact that
both pairs of functions $(a,b)$ and $(a',b')$ determine the same distributions.
$\Box$

\smallskip
{\bf Theorem 4.}
{\it The ODE {\rm (\ref{du2})} of hyperbolic type is equivalent to
two pairs of functions  $(a,b)$, $(f,g)$ on  $M$ in general
position.}


\section{The main result}

Now we return to the ODEs $y''=\lambda_i(x,y,y')$ and vector
fields $X_i$, $i=1,2$. Choose two pairs of independent integrals
$(a, b)$ and $(f,g)$ of  $X_1$ and $X_2$ respectively. Let
$g=h(a,b,f)$. Let
$$
a, ~b, ~c=-\dfrac{h_{a}}{h_{b}}
$$
be new coordinates in
$\mathbb{R}^3$ and let  $G_2(a, b,c)=- {{\rm Flex} (h)}/{h_{b}^3}$.
Then $X_1$ and $X_2$ will have the form
\begin{equation}
\label{X1X2}
X_1=\dfrac{\partial}{\partial c}
~~\hbox{and}~~
X_2=\dfrac{\partial}{\partial a} +
c\dfrac{\partial}{\partial b}+G_2(a,b,c)
\dfrac{\partial}{\partial c},
\end{equation}
respectively.
We will say that  $(a, b, c)$ are  {\it
canonical coordinates} and that the ODE
\begin{equation}
\label{eqG}
 b''=G_2( a, b, b'), \quad b=b(a),
\end{equation}
determined by $X_2$ is {\it associated} with the ODE~(\ref{du2}).

Let $M_3$ and $\widetilde M_3$ be two 3-dimensional contact
manifolds.
Consider the contact transformation
$\varphi:M_3\to \widetilde M_3$ which maps
the ODE (\ref{du2}) into another equation of the same form.
Let it map the vector field $X_1$
to the vector field $\widetilde X_1=\varphi^*X_1$.
We prove that $\varphi$ acts of
associated equations as a point transformation.

Let $M_2=M_3/X_1$ be the factor of $M_3$ by the trajectories of
$X_1$   (locally it can be done).
Let  $\pi:M_3\to M_2$ be the canonical projection.
Similarly, let $\widetilde{\pi}:\widetilde M_3\to
\widetilde M_2=\widetilde M_3/\widetilde{X}_1$.
The modelling bundle for  $\pi$ is the bundle
$\pi_0:J^1N\to J^0N$ where $N$ is some 1-dimensional manifold.

Let $\varphi_0:J^1N\to J^1\widetilde N$ be the restriction on
$\varphi$ to $J^1N$.
The diffeomorphism $\varphi$ maps the integrals
$a$ and $b$ of the distribution determined by $X_1$
to some integrals $\widetilde a$ and $\widetilde b$
of the distribution determined by $\widetilde X_1$,
i.e., to functions of $a$ and~$b$.
Thus, it induces the diffeomorphism
$\psi_0:J^0N\to J^0\widetilde N$, which closes the diagram
$$
\xymatrix@R=15mm{%
    {J^1N} \ar[rr]^-{\varphi_0}\ar[d]^{\pi_0} & &
         {J^1\widetilde N}\ar[d]_{\widetilde\pi_0}  \\
    {J^0N} \ar[rr]^-{\psi_0}  &
     & {J^0\widetilde N.}  \\
}
$$
Since one can choose the functions $a$ and $b$ as coordinates on
$J^0N$, one can treat  $\psi_0$ as a point transformation.
Conversely, any diffeomorphism $\psi_0:J^0N\to J^0\widetilde N$
lifts to the contact diffeomorphism $\varphi_0:J^1N\to
J^1\widetilde N$. Hence, the canonical coordinate system is
defined up to a change of $(a, b)$ to the pair of independent
functions $(\widetilde a(a,b), \widetilde b(a,b))$.

\smallskip
{\bf Proposition 2.}
{\it
The point equivalence class of the associated
equation {\rm(\ref{eqG})} does not depend on the choice of
the integrals $(a, b)$, $(f, g)$.}

{\bf Proof.}
The choice of another integrals $\widetilde a$, $\widetilde b$ of $X_1$
instead of  $a$, $b$ induces the change of canonical coordinates,
thus does not change the equivalence class of the ODE {\rm(\ref{eqG})}.

Let $F(a,b,f,g)=0$ denote some functional dependency
between these integrals. Any other
dependency (in particular,
$g-h(a,b,f)=0$) has the form
$$
\widetilde F(a,b,f,g)=\Phi(F(a,b,f,g))=0.
$$
Clearly, the ratio ${\widetilde F_a}/{\widetilde F_b}$ does not depend on $\Phi$, in
particular,
$$
\frac{F_a}{F_b}=\frac{h_a}{h_b}.
$$
The choice of another integrals $\widetilde f=\phi(f,g)$ и
$\widetilde g=\psi(f,g)$ also does not change the ratio
${F_a}/{F_b}$,
hence, does not change the canonical coordinate system.
$\Box$
\smallskip

In the construction above we can change the order of
$X_1$ and $X_2$.
Thus we obtain another ODE
$b''=G_1(a, b, b')$, associated
with {\rm (\ref{du2})}.
We say that these ODEs are {\it dual} to each other
(see~\cite{Ar}).

From above it follows the

\smallskip
{\bf Theorem 5.} {\it Let $\mathcal{E}$ and
$\widetilde{\mathcal{E}}$ be two ODEs {\rm (\ref{du2})} of
hyperbolic type and  let $([\mathcal{E}^a_1], [\mathcal{E}^a_2])$
and $([\widetilde{\mathcal{E}}^a_1],
[\widetilde{\mathcal{E}}^a_2])$ be the pairs of point equivalence
classes of their associated ODEs. The ODEs $\mathcal{E}$ and
$\widetilde{\mathcal{E}}$ are contact equivalent if and only if
one of the classes $([\mathcal{E}^a_1], [\mathcal{E}^a_2])$
coincides with one of the classes $([\widetilde{\mathcal{E}}^a_1],
[\widetilde{\mathcal{E}}^a_2])$.}

\smallskip
Let us present some examples.

{\bf Example 1.}
Consider the equation
\begin{equation}
\label{ex1}
y''^2-y''=0.
\end{equation}
Here $\lambda_1=1$, $\lambda_2=0$.
Let
$$
a=p-x,~~ b=\frac12x^2-px+y, ~~f=p,
~~
g=y-px.
$$
Then $g=b-\frac12(a-f)^2$, which yields $G_2=1$.
On the other hand,  $b=g+\frac12(a-f)^2$, which yields
$G_1=-1$.
Thus, the pair of ODEs associated with
(\ref{ex1}) are
$$b''=1
~~\hbox{ and }~~
b''=-1.
$$
Clearly, they are point equivalent.

{\bf Example 2.}
Consider the equation
$$
(y''-y)(y''-y')=0.
$$
Here we have $\lambda_1=y$, $\lambda_2=p$.
Let
$$
a=p^2-x^2, ~~b=x-\ln(p+y),~~
f=y-p, ~~g=pe^{-x}.
$$
The pair of associated ODEs is
$$
b''=-\dfrac{b'(2+ab')}a  ~~~~\hbox{and}~~~~
b''=\dfrac{b'}{a}.
$$
Both of them are point equivalent to the ODE $b''=0$.

{\bf Example 3.}
Consider the equation
$$
(y''+x)(y''-1)=0,
$$
for which $\lambda_1=-x$, $\lambda_2=1$.
Let
$$
a=\frac12x^2+p, ~~b=\frac13x^3-px+y,
~~
f=p-x, ~~g=\frac12x^2-px+y.
$$
The pair of associated ODEs is
$$
b''=\dfrac{1}{b'-1}  ~~~~\hbox{and}~~~~
b''=\dfrac{-1-4b'{}^2+2\sqrt{b'{}^2+4b'{}^4}}{(1+8b'{}^2
-4\sqrt{b'{}^2+4b'{}^4})^{3/2}}.
$$
Both of them admit the same Lie algebra, hence they are point
equivalent.

{\bf Example 4.}
We give the example of two pairs of functions in general position
for which associated ODEs are not point equivalent to each other.
Let functions $(a,b)$ and $(f,g)$ satisfy the conditions
$$
a^2 - 2a b + f^2 + 2 b f - g^2 + b^2 = 1,
\quad g\ge0, \quad b-a+f\ge0.
$$
Then $g=h(a,b,f)=\sqrt{a^2 - 2a b + f^2 + 2 b f  + b^2 - 1}$,
and the associated equation is
$$
b''=\dfrac{b'}{a - b} - \dfrac{2b'{}^2}{a - b} +
\dfrac{b'{}^3}{a - b}.
$$
It is point equivalent to the equation $b''=0$.

On the other hand, $b=h(f,g,a)= a - f + \sqrt{1 - 2 a f + g^2}$,
and the second associated equation is not cubic
in the first derivative, hence it is not point equivalent to the
first one.

\smallskip
In conclusion we give the formulas needed to calculate the
invariants used in Tresse-Kruglikov theorem.
We denote $\dfrac{\partial}{\partial a}$,
$\dfrac{\partial}{\partial  b}$ and
$\dfrac{\partial}{\partial  c}$
by $\nabla_{a}$, $\nabla_{b}$ and $\nabla_{c}$.
In the initial coordinates
$(x,y,p)$ they are of the form

$$
\begin{array}{l}
\nabla_{a}=\dfrac1{\Delta}
\Bigl((b_{y}f_{p} - b_pf_y)
\dfrac{\partial}{\partial x}
+ (b_{p}f_x-b_xf_p)
\dfrac{\partial}{\partial y}
+(b_xf_y-f_xb_y)
\dfrac{\partial}{\partial p}\Bigr),\\[10pt]
\nabla_{b}=\dfrac1{\Delta}\Bigl(
( a_p f_y  - a_y f_p)
\dfrac{\partial}{\partial x}
+
 (a_xf_p -a_pf_x )
\dfrac{\partial}{\partial y}
+(a_y f_x - a_x f_y)
\dfrac{\partial}{\partial p}\Bigr),\\[10pt]
\nabla_{c}=\dfrac1{\Delta}\Bigl(
( a_y b_p - a_p b_y )
\dfrac{\partial}{\partial x}
+
(  a_pb_x - a_{x} b_p)
\dfrac{\partial}{\partial y}
+
(a_xb_y - a_y b_x)
\dfrac{\partial}{\partial p}\Bigr),
\end{array}
$$
where
$$
\Delta=
a_pb_xf_y-a_pb_yf_x-a_yp_xf_p+a_yb_pf_x+a_xb_yf_p-a_xb_pf_y,
$$
and the partial derivatives  of $f$ in $x$, $y$ and $p$
are understood as derivatives of a composite function.

The invariant $H$ is
\begin{multline*}
\t H=\nabla_c^2\nabla_a^2F  - 2\cdot\dfrac{h_a}{h_b}\cdot\nabla_c^2\nabla_a\nabla_b F
+\left(\dfrac{h_a}{h_b}\right)^2\cdot \nabla_c^2\nabla_b^2 F  -4\nabla_c\nabla_a\nabla_b  F
+ 4\cdot\dfrac{h_a}{h_b}\cdot \nabla_c\nabla_b^2F +\\
+6 \nabla_b^2 F
 + F\left(2\nabla_c^3\nabla_a F-2\cdot\dfrac{h_a}{h_b}\cdot \nabla_c^3\nabla_b F
-  3\nabla_c^2\nabla_bF\right)  -\nabla_cF\biggl(\nabla_c^2\nabla_aF-\dfrac{h_a}{h_b}\cdot
\nabla_c^2\nabla_bF - \\
-
4\nabla_c\nabla_bF\biggr) +  \nabla_c^3F \left(\nabla_aF-
\dfrac{h_a}{h_b}\cdot \nabla_bF\right) -
 3\nabla_c^2F \cdot\nabla_bF+ F^2\cdot\nabla_c^4 F.
\end{multline*}

Invariant differential operators are
$$
\begin{array}{l}
\Delta_c=\widetilde{\cal D}_c +
(r-s)\dfrac{\nabla_c^5F}{5\nabla_c^4F},
\\[7pt]
\Delta_a = \widetilde{\cal D}_a -\dfrac{h_a}{h_b}\cdot
\widetilde{\cal D}_b
+\,F\widetilde{\cal D}_c + (3r+2s)\nabla_cF +
\\[4pt]
\qquad +\,(2r+s)\Biggl(\dfrac{F\cdot
\nabla_c^5F+
\nabla_c^4\nabla_aF -\, \dfrac{h_a}{h_b}\cdot\nabla_c^4\nabla_bF}{\nabla_c^4 F}
\Biggr),
\\[7pt]
\Delta_b=  \dfrac{\nabla_c^5F}{5\nabla_c^4F}\cdot\widetilde{\cal
D}_a
+ \Biggl(1-\dfrac{h_a}{h_b}\cdot\dfrac{\nabla_c^5F}{5\nabla_c^4F}
\Biggr)\widetilde{\cal D}_b\,+
\\[4pt]
\qquad+\,
\Biggl(2\nabla_cF+ \dfrac{6F\cdot
\nabla_c^5F + 5\nabla_c^4\nabla_aF
-5\cdot\dfrac{h_a}{h_b}\cdot\nabla_c^4\nabla_bF}{\nabla_c^4 F}
\Biggr)
+(3r+2s) \dfrac{\nabla_c^2 F}{8} +
\\[7pt]
\qquad+\,(r+2s)\dfrac{\nabla_c^4\nabla_bF}{4\nabla_c^4F} +
(19r+6s)\dfrac{\nabla_c^5F\cdot\nabla_cF}{10\nabla_c^4F}
+ (21r+6s)\dfrac{\nabla_c^5F(F\cdot \nabla_c^5F +
\nabla_c^4\nabla_bF)}{20\nabla_c^4F\cdot \nabla_c^4F}
 ,
\end{array}
$$
where  $\widetilde{\cal D}_a$, $\widetilde{\cal D}_b$,
$\widetilde{\cal D}_c$ --- are the total derivatives corresponding to
the operators $\nabla_{a}$, $\nabla_{b}$,
$\nabla_{c}$.

{\bf Acknowledgement.}
Author wishes to express his deep gratitude to V.V.\,Lychagin for
the formulation of the problem and  valuable comments.

\end{document}